\documentclass[12pt]{article}
\usepackage{amsmath}
\usepackage{amsfonts}
\usepackage{amssymb}
\usepackage{theorem}

\title{Sofic metric groups and continuous logic}
\author{A.Ivanov 
\thanks{
The research is supported by Polish National Science Centre grant DEC2011/01/B/ST1/01406
} 
}

\newtheorem{theorem}{Theorem}[section]
\newtheorem{proposition}[theorem]{Proposition}

\newtheorem{definition}[theorem]{Definition}

\newtheorem{remark}[theorem]{Remark}

\begin{document} 
\topmargin = 12pt
\textheight = 630pt 
\footskip = 39pt 

\maketitle

\begin{quote}
{\bf Abstract}
We define sofic, weakly sofic, linear sofic and hyperlinear 
metric groups and discuss some issues involving axiomatizability 
of these classes in continuous logic. \\ 
{\bf 2010 Mathematics Subject Classification}: 20E26; 03C20.\\ 
{\bf Keywords}: Sofic groups, continuous logic. 
\end{quote}

\section{Introduction}
\label{intro}

This paper concerns several questions involving 
axiomatizability of the class 
of sofic groups in continuous logic. 
This issue has been already discussed in our 
paper \cite{sasza} (see the final part of Introduction).  
Now we extend this discussion in order to answer 
several questions addressed to the author 
after posting \cite{sasza}. \parskip0pt

Let us consider the class $\mathcal{G}$ of all continuous 
structures which are metric groups $(G,d)$ with 
bi-invariant metrics $d\le 1$. 
We defined in \cite{sasza} the subclass 
$\mathcal{G}_{sof} \subset \mathcal{G}$ 
of all {\em sofic metric groups} as all closed metric 
subgroups of metric ultraproducts of 
finite symmetric groups with Hamming metrics. 
In order to approach to the Gromov's question of 
soficity of all countable groups (see \cite{CL}, \cite{pestov}) 
we asked if the equality 
$\mathcal{G}_{sof} = \mathcal{G}$ 
holds.  
We emphasize that in this question groups are 
considered together with metrics. 
Thus this equality is a stronger version of 
the problem of soficity of all countable groups. 
We hoped that the corresponding counterexamples 
could be interesting with respect to the original question.  
We will see below that in fact there are 
easy examples of (finite) metric groups 
which do not belong to  $\mathcal{G}_{sof}$. 
At the moment we do not know how useful they are. 

On the other hand let us consider 
the class $\mathcal{G}_{w.sof}$ 
of weakly sofic continuous metric groups, 
i.e. continuous metric groups $(G,d)$ which embed  
into a metric ultraproduct of finite metric groups with 
invariant length functions bounded by 1 \cite{GR}. 
We now see that it properly extends 
$\mathcal{G}_{sof}$. 
Does this class coincide with $\mathcal{G}$?

Note that the positive answer to this 
question implies that any group is weakly sofic 
in  the standard sense of \cite{GR}. 
On the other hand if there is a non-weakly sofic 
group $G$ in the standard sense, then for every 
metric $d$ the metric group $(G,d)$ does not satisfy 
a sentence of continuous logic which holds 
in all weakly sofic metric groups. 
This follows from the observation that 
$\mathcal{G}_{w.sof}$ is the minimal subclass of 
$\mathcal{G}$ which is axiomatizable in continuous 
logic and which contains 
all finite metric groups with metrics $\le 1$ 
(see the corresponding discussion in Section 2). 

It is interesting that to get a positive answer to 
the questions above it suffices to prove that 
all finitely generated free groups 
(i.e. sofic groups) equipped with bi-invariant discrete 
metrics $\le 1$  belong to  
$\mathcal{G}_{w.sof}$.  
This follows from \cite{Doucha} (see Section 3.2 there 
and Remark \ref{doucha} below). 

In our paper we discuss these and some other 
questions of this type, for example in the cases 
of linear sofic groups and hyperlinear  groups. 
We pay a special attention to discrete members of these classes. 
All our results will be formulated in Section 2. 
It can be considered as a continuation of this introduction. 

In the rest of the introduction we briefly remind 
the reader some preliminaries of continuous logic. 
Although this material is not necessary for our main results 
it serves as a platform of the paper. 
In particular Theorem \ref{axiom} and its corollary 
given in Proposition  \ref{sup-ax} tell us which questions 
are basic in our investigations.

\bigskip

\paragraph{Continuous structures.} 

We fix a countable continuous signature 
$$
L=\{ d,R_1 ,...,R_k ,..., F_1 ,..., F_l ,...\}. 
$$  
Let us recall that a {\em metric $L$-structure} 
is a complete metric space $(M,d)$ with $d$ bounded by 1, 
along with a family of uniformly continuous operations on $M$ 
and a family of predicates $R_i$, i.e. uniformly continuous maps 
from appropriate $M^{k_i}$ to $[0,1]$.   
It is usually assumed that to a predicate symbol $R_i$ 
a continuity modulus $\gamma_i$ is assigned so that when 
$d(x_j ,x'_j ) <\gamma_i (\varepsilon )$ with $1\le j\le k_i$ 
the corresponding predicate of $M$ satisfies 
$$ 
|R_i (x_1 ,...,x_j ,...,x_{k_i}) - R_i (x_1 ,...,x'_j ,...,x_{k_i})| < \varepsilon . 
$$ 
In this paper we usually assume that $\gamma_i$ coincides with $id$. 
In this case we do not mention the appropriate modulus. 
We also fix continuity moduli for functional symbols. 
Note that each countable structure can be considered 
as a complete metric structure with the discrete $\{ 0,1\}$-metric. 

By completeness continuous substructures of a continuous structure are always closed subsets. 

Atomic formulas are the expressions of the form $R_i (t_1 ,...,t_r )$, 
$d(t_1 ,t_2 )$, where $t_i$ are terms (built from functional $L$-symbols). 
In metric structures they can take any value from $[0,1]$.   
{\em Statements} concerning metric structures are usually 
formulated in the form 
$$
\phi = 0 
$$ 
(called an $L$-{\em condition}), where $\phi$ is a {\em formula}, 
i.e. an expression built from 
0,1 and atomic formulas by applications of the following functions: 
$$ 
x/2  \mbox{ , } x\dot- y= max (x-y,0) \mbox{ , } min(x ,y )  \mbox{ , } max(x ,y )
\mbox{ , } |x-y| \mbox{ , } 
$$ 
$$ 
\neg (x) =1-x \mbox{ , } x\dot+ y= min(x+y, 1) \mbox{ , } sup_x \mbox{ and } inf_x . 
$$ 
Note that all diadic numbers $\le 1$ are formulas. 
A {\em theory} is a set of $L$-conditions without free variables 
(here $sup_x$ and $inf_x$ play the role of quantifiers). 
   
It is worth noting that any formula is 
a $\gamma$-uniformly continuous function 
from the appropriate power of $M$ to $[0,1]$, 
where $\gamma$ is the minimum of continuity 
moduli of $L$-symbols appearing in the formula. 

The condition that the metric is bounded by $1$ 
is not necessary. 
It is often assumed that $d$ is bounded 
by some rational number $d_0$. 
In this case the (truncated) functions 
above are appropriately modified.

Given a family $M_i$ of structures with metrics $d_i$, $i\in I$, 
and an ultrafilter $D$ on $I$, 
the {\em metric ultraproduct}  
$\prod_I  M_i /D$ is defined as follows. 
For  $( u_i )_I$ and $( v_i )_I \in \prod_I  M_i$ 
define the distance  
$$ 
d(( u_i )_I , ( v_i )_I ) = 
lim_{i \rightarrow D} d_i (u_i ,v_i ) , 
$$ 
i.e. by the rule that for any $\varepsilon_1 <\varepsilon_2$ 
from $[0,1]$ the distance between 
$(u_i )_I$ and $(v_i )_I$ is in the interval 
$(\varepsilon_1 , \varepsilon_2 )$  if and only if 
the set $\{ i : d_i (u_i ,v_i ) \in (\varepsilon_1 , \varepsilon_2 )\}$ 
belongs to the ultrafilter $D$. 
This is well-defined by compactness of $[0,r]$. 
The structure $\prod_I  M_i /D$ consists of 
classes of the relation 
$d(( x_i )_I , ( y_i )_I ) = 0$. 
Its operations and predicates are defined in 
the natural way.

\paragraph{Axiomatisability in continuous logic} 

When one considers classes axiomatizable 
in continuous logic it is obviously assumed 
that all operations and predicates are uniformly 
continuous with respect to some fixed continuity moduli. 
Suppose $\mathcal{C}$ is a class of metric $L$-structures. 
Let $Th^c (\mathcal{C})$ be the set of all closed $L$-conditions 
which hold in all structures of $\mathcal{C}$. 
It is proved in \cite{BYBHU} (Proposition 5.14 and Remark 5.15) 
that every model of $Th^c (\mathcal{C})$ is elementary equivalent 
to some ultraproduct of structures from $\mathcal{C}$. 
Moreover by Proposition 5.15 of \cite{BYBHU} we have the following statement. 
\begin{quote} 
The class $\mathcal{C}$ is axiomatizable in continuous logic 
if an only if it is closed under metric isomorphisms and 
ultraproducts and its complement is closed under ultrapowers. 
\end{quote} 
Let $Th^c_{\sup} (\mathcal{C})$ be the set of all closed 
$L$-conditions of the form 
$$ 
sup_{x_1} sup_{x_2} ... sup_{x_n} \varphi =0 
\mbox{ ( $\varphi$ does not contain $inf_{x_i}$, $sup_{x_i}$ ), } 
$$
which hold in all structures of $\mathcal{C}$.
Some standard arguments also give the following theorem. 

\begin{theorem} \label{axiom} 
(1) The class $\mathcal{C}$ is axiomatizable in continuous logic 
if an only if it is closed under metric isomorphisms, ultraproducts and 
taking elementary submodels. \parskip0pt  

(2) The class $\mathcal{C}$ is axiomatizable in continuous logic 
by $Th^c_{sup} (\mathcal{C})$ if an only if it is closed under 
metric isomorphisms, ultraproducts and taking substructures. 
\end{theorem} 

{\em Proof.} 
Statements (1) and (2) are similar. 
In both cases necessity is obvious. 
Let us consider the sufficiency direction of statement (2). 
Let a continuous structure $M$ 
satisfy $Th^c_{sup} (\mathcal{C})$. 
The atomic diagram of $M$ consists of all 
equalities of the form $d(t_1 ,t_2 ) =r$
and $R(\bar{t}) = r$ with parameters from $M$. 
Fixing a small rational $\varepsilon >0$ 
we replace these equalities 
by the corresponding inequalities
$$ 
|d(t_1 ,t_2) - r' |\le \varepsilon /2 \mbox{ , } 
|R(\bar{t}) - r' | \le \varepsilon /2 , 
$$
where $r'$ is a rational number 
with $| r-r' | < \varepsilon /2$.  
Let $D_{\varepsilon}(M)$ be the set of all 
statements of this form.  
Let $I$ be the collection of all finite 
subsets of the union 
$\bigcup \{ D_{\varepsilon}(M) : \varepsilon \in \mathbb{Q}^+ \}$. 

Let $i \in I$ depend on parameters $\bar{a}$. 
We may assume that $i\subset D_{\varepsilon}(M)$ 
for some fixed $\varepsilon$ (choosing the minimal one in $i$). 
Let $\bar{x}$ correspond to parameters $\bar{a}$. 
Note that the statemet 
$$ 
\sup_{\bar{x}} max ( \{ |d(t_1 ,t_2) - r' | : \mbox{ an inequality of the form } 
|d(t_1 ,t_2) - r' |\le \varepsilon' /2 
$$
$$
\mbox{ appears in } i (\bar{x})  \} \cup  
\{ |R(\bar{t}) - r' | : \mbox{ an inequality of the form } 
$$
$$
|R(\bar{t}) - r' | \le \varepsilon' /2 \mbox{ appears in }  i(\bar{x}) \})  
\ge \varepsilon /2 
$$
is not realized by $\bar{a}$ in $M$. 
Since $M \models Th^c_{sup} (\mathcal{C})$, 
there are $M_i \in \mathcal{C}$ and $\bar{a}_i \in M_i$ 
so that $i (\bar{x})$ is realized in $M_i$ by  $\bar{a}_i$. 

Let $J_i = \{ j\in I : i\subseteq j \}$ and 
let $\mathcal{U}$ be an ultrafilter over $I$ containing 
all $J_i$. 
The structure $M$ can be embedded into  
$$
\prod_{I} M_i /\mathcal{U}  
$$ 
as follows. 
If $c$ appears in $i$, then we map $c$ to the 
corresponding element of the tuple $\bar{a}_i$. 
It is easy to see that the atomic diagram of $M$
is satisfied by the induced mapping 
into $\prod_{I} M_i /\mathcal{U}$. 
$\Box$

\paragraph{Metric groups} 
Let us consider continuous metric structures 
which are groups. 
Below we always assume that all groups which are 
continuous structures are considered with respect 
to bi-invariant metrics. 
We take this assumption because any group which is 
a continuous structure has such a metric. 
See \cite{sasza2} for a discussion 
concerning  this observation. 

Since the algebraic approach from \cite{G}, 
\cite{pestov} and \cite{thom} has become dominating 
in the literature, we now consider metric groups 
in slightly different terms. 

Let $G$ be a group. 
A function $l:G\rightarrow [0,\infty)$ is called 
a {\em pseudo length function} if \\ 
(i) $l(1) = 0 $; \\ 
(ii) $l(g)=l(g^{-1})$; \\
(iii) $l(gh)\le l(g) + l(h)$. \\ 
A {\em length function} is a pseudo length 
function satisfying 
\begin{quote} 
$(i')$ $l(g) = 0 $ if and only if $g=1$, where $g\in G$.  
\end{quote} 
A pseudo length function is {\em invariant} if 
$l(h^{-1} gh) = l(g)$ for all $g,h\in G$. 
In this case it defines an invariant 
pseudometric $d$ by $l(gh^{-1})$. 
It becomes a metric if $l$ is a length function. 
In this case we say that $G$ is a {\em normed group}. 
We consider normed groups with bounded norms as 
a notion which is equivalent to  metric groups 
(an easy exercise). 
It is worth noting that any unbounded 
bi-invariant norm $l$ can be replaced by 
the norm $h \rightarrow \frac{l(h)}{1+l(h)}$ 
which defines the same topology with $l$.

Metric ultraproducts of normed groups of 
bounded diameter, say $r$, are defined as follows.  
Let $(G_i ,l_i )$, $i\in I$, be a family of groups 
equipped with invariant length functions 
and let $\Phi$ be an ultrafilter on $I$. 
Then 
$$ 
N=\{ (x_i )_{i\in I}\in \prod_I G_i : 
lim_{i\rightarrow  \Phi} l(x_i ) = 0 \} 
$$ 
is a normal subgroup of $\prod_I G_i$.   
The {\em metric ultraproduct} $\prod_I (G_i ,l_i )/\Phi$ 
is defined to be $(\prod_I G_i )/N$ where 
the length function is defined by 
$$ 
l(xN) = lim_{i\rightarrow \Phi} l_i (x_i ) . 
$$ 
This definition corresponds to Section 2.4 from \cite{pestov} 
and  to the  definition of 
metric ultraproducts in continuous model theory
which was given above.

Since axiomatizable closures of classes of discrete structures 
are central in our paper, we make two easy remarks concerning 
this situation. 

Firstly note that the class of all groups with the $\{ 0,1\}$-metric 
is axiomatizable both in the discrete and in the continuous logic
\footnote{in this case axiomatizability in continuous logic is equivalent to axiomatizability in first-order logic} 
.  
If we consider the class of all abstract groups under all discrete 
metrics then it becomes non-axiomatizable.   
For example it is not closed under metric ultraproducts.  

On the other hand it may also happen that when 
we extend an axiomatizable class of groups with the $\{ 0,1\}$-metric 
by (abstract) structures from this class with all possible 
(not only possible discrete) metrics we lose axiomatizability. 
A nice example of this situation is the class of non-abelian groups 
with $[0,1]$-metrics. 
For example there is a sequence of non-abelian groups $G_n \le Sym (2^n +3 )$ 
with $G_n \cong \mathbb{Z}(2)^n \times S_3$ so that their metric unltraproduct 
with respect to Hamming metrics is abelian (an easy exercise).

\section{Axiomatizability in continuous logic and  sofic groups } 
\label{SecAx}

Metric ultraproducts of finite normed groups are  
deserved a particular attention in group theory.  
This is mainly motivated by investigations of 
{\em sofic groups}.  
We remind the reader that a group $G$ is called 
{\em sofic} if $G$ embeds into 
a metric ultraproduct of finite symmetric groups 
with the Hamming distance $d_{H}$, \cite{pestov}. 
We remind the reader that 
$$
d_H (g,h) = 1 - \frac{|Fix (g^{-1}h)|}{n} \mbox{ for } 
g,h \in S_n . 
$$ 
A group $G$ is called {\em hyperlinear} if $G$ embeds into 
a metric ultraproduct of finite-dimensional unitary  groups  $U(n)$ 
with the normalized Hilbert-Schmidt metric $d_{HS}$ 
(i.e. the standard $l^2$ distance between matrices), \cite{pestov}. 
It is an open question whether these classes are the same 
and whether any countable group is sofic/hyperlinear.

In the case of groups it is convinient to use the following notion 
of approximation, see  
\cite{thom} and \cite{G} (Definition 3). 

\begin{definition} \label{DefApp} 
{\em 
Let $\mathcal{K}$ be a class of groups and $\mathcal{L}$ 
be a class of invariant length functions on groups 
from $\mathcal{K}$ 
({\em we will always assume that they are bounded by some fixed number} 
$r$). 
We say that a group $G$ is $(\mathcal{K},\mathcal{L})$-approximable if 
there is a function $\alpha : G\rightarrow [0,r]$ with 
$$ 
\alpha(g) = 0 \Leftrightarrow g=1 , 
$$ 
so that for any finite $F\subset G$ and $\varepsilon >0$ there is 
$(H,l) \in \mathcal{K}$, $l\in \mathcal{L}$, and a function $\gamma : F \rightarrow H$ 
so that 
$$
l(\gamma (1)) <\varepsilon \mbox{ , } \forall g\in F (l(\gamma (g))\ge \alpha (g)) \mbox{ , and }  
$$
$$ 
l (\gamma (gh) (\gamma (g) \gamma(h) )^{-1} ) < \varepsilon \mbox{ , for any } g,h, gh \in F . 
$$ 
}
\end{definition} 

It is known that a ghroup $G$ is $(\mathcal{K},\mathcal{L})$-approximable 
if and only if it embeds into a metric ultraproduct of groups 
from $\mathcal{K}$ with norms from $\mathcal{L}$ (\cite{thom} and \cite{G}). 
Moreover in the case of sofic and hyperlinear groups the function $\alpha$ 
can be taken constant on $G\setminus \{ 1 \}$ with the value 
equal to any real number strictly between $0$ and $1$ 
(between $0$ and $2$ in the hyperlinear case).  
 
A group $G$ is called {\em weakly sofic} if $G$ embeds 
into a metric ultraproduct of finite  groups with 
invariant length functions bounded by 1 \cite{GR}. 
It is not known if this class coincides with the former ones.

When we consider {\em metric groups} this definition of approximation 
should be modified as follows. 

\begin{definition} \label{DefMetrApp} 
{\em 
Let $\mathcal{K}$ be a class of groups and $\mathcal{L}$ 
be a class of invariant length functions on groups 
from $\mathcal{K}$. 
We say that a metric group $(G,d)$ is  
$(\mathcal{K},\mathcal{L})$-approximable if 
for any finite $F\subset G$ and $\delta >0$ there is 
$(H,l) \in \mathcal{K}$, $l\in \mathcal{L}$, and a function 
$\gamma : F \rightarrow H$ so that 
$$  
(\forall g,h, gh \in F )(
l (\gamma (gh) (\gamma (g) \gamma(h) )^{-1} ) < \delta )
\mbox{ , if } 1 \in F \mbox{ , then } l(\gamma (1) ) < \delta \mbox{ , } 
$$
$$ 
 \mbox{ and }( \forall g,h\in F)(| d(g,h) - l(\gamma (g) \gamma (h)^{-1})| <\delta   . 
$$
} 
\end{definition} 

It is easy to see that a metric group $(G,d)$ with $d\le 1$ is  
$(\mathcal{K},\mathcal{L})$-approximable if and only if 
it isometrically embeds into a metric ultraproduct of $\mathcal{L}$-normed 
groups from $\mathcal{K}$.

Let us recall that $\mathcal{G}_{sof}$ denotes the class of 
metric groups of diameter 1 with bi-invariant metrics  
so that the groups are embeddable as closed subgroups  
via isometric morphisms into a metric ultraproduct of 
finite symmetric groups with Hamming metrics.  
In particular $\mathcal{G}_{sof}$ consists of 
continuous metric structures. 
We call it the class of {\em sofic metric groups}.

By $\mathcal{G}_{w.sof}$ and  $\mathcal{G}_{hyplin}$ 
we denote the classes of continuous structures which are 
weakly sofic metric groups and hyperlinear metric groups 
(which are considered under the metric induced by 
$\frac{1}{2}d_{HS}$ in the corresponding ultraproducts of $U(n)$'s). 
In \cite{ArP2} Arzhantseva and P\"{a}unescu introduced 
{\em linear sofic groups} as ones which are approximated by 
all $GL_n (\mathbb{C})$ considered together with 
the metric $\rho (a,b) = n^{-1}\mbox{ rk}(a-b)$. 
Applying Definition \ref{DefMetrApp} we introduce 
the class $\mathcal{G}_{l.sof}\subset \mathcal{G}$ of 
linear sofic metric groups.   

\begin{proposition} \label{sup-ax} 
The class of sofic  
(resp. linear sofic, hyperlinear and weakly sofic) 
metric groups is $sup$-axiomatizable, 
i.e. by its theory $Th^c_{sup}$.  
\end{proposition} 

{\em Proof.} 
The class of metric ultraproducts of $(S_n ,d_H )$, $n\in \mathbb{N}$, 
is closed under iteration of metric ultraproducts. 
Thus the case of $\mathcal{G}_{sof}$ follows from Theorem \ref{axiom}. 
The remained cases are similar. 
$\Box$ 
\bigskip 
 
Let us consider relationship among the classes 
of the following collection 
$$
\{ \mathcal{G} ,  \mathcal{G}_{sof} , \mathcal{G}_{w.sof} ,  \mathcal{G}_{hyplin},   
\mathcal{G}_{l.sof} \} 
$$ 
It is clear that 
$\mathcal{G}_{sof} \subseteq  \mathcal{G}_{w.sof} \subseteq \mathcal{G}$.  
Moreover, by the proof of Theorem 8.2 of \cite{ArP2} 
$\mathcal{G}_{l.sof} \subseteq \mathcal{G}_{w.sof}$. 
Although the arguments from \cite{pestov} and \cite{ArP2} 
that sofic groups are hyperlinear and linear sofic, 
do not allow the corresponding approximations of Hamming metrics, 
it still looks likely that $\mathcal{G}_{sof}$ 
is contained in the other classes. 
We will see below that there is an example from 
$\mathcal{G}_{w.sof}\setminus (\mathcal{G}_{sof} \cup  \mathcal{G}_{hyplin} \cup\mathcal{G}_{l.sof})$. 
Thus $\mathcal{G}_{w.sof}\not\subset \mathcal{G}_{sof}$. 
In order to approach the other inclusions 
we investigate how discrete structures of these classes 
represent them. 

It is folklore that any abstract sofic group can be 
embedded into a metric ultraproduct of finite symmetric 
groups as a discrete subgroup 
(see the proof of Theorem 3.5 of \cite{pestov}). 
This means that the class of all abstract sofic groups consists 
of all groups which can be viewed as discrete structures of 
the class $\mathcal{G}_{sof}$ (with the $\{ 0,1\}$-metric).    
On the other hand any weakly sofic group 
(or weakly sofic metric group)
embeds into a metric ultaproduct 
of discrete (even finite) metric groups. 
The following theorem is related to these observations 
in the cases of all groups and hyperlinear metric groups. 

\begin{theorem} \label{main} 
Let $(G,d)$ be a bi-invariant metric group so that  
$d\le 1$ and $(G,d)$ is a continuous structure. 
Then the following statements hold. 

1. $(G,d)$ is a closed subgroup of a metric 
ultraproduct of discrete bi-invriant metric groups. 
Moreover if $G$ is weakly sofic in the standard 
sense of \cite{GR}, then the groups occurring 
in the corresponding ultaproduct are 
weakly sofic too.

2.  If $(G,d)$ is hyperlinear 
then  $(G,d)$ is a closed subgroup of a metric ultraproduct of 
discrete bi-invriant metric groups which are hyperlinear. 
\end{theorem} 

The theorem will be proved in the next section. 
It shows that if two 
classes $\mathcal{K}_1$ , $\mathcal{K}_2$ from the collection 
$$
\{ \mathcal{G} ,  \mathcal{G}_{sof} , \mathcal{G}_{w.sof} ,  \mathcal{G}_{hyplin},   
\mathcal{G}_{l.sof} \} 
$$ 
have the same countable discrete structures, then  
$\mathcal{K}_1 = \mathcal{K}_2$. 
Indeed  if $\mathcal{K}_1$ and $\mathcal{K}_2$ have 
the same discrete structures then they are generated as 
axiomatizable classes by the same set of structures 
(in the linear sofic case note that $(GL_n (\mathbb{C}), \rho )$ 
is a discrete structure). 
Thus they are the same. 
On the other hand any discrete metric structure of a finite 
language is naturally embedded into a metric ultraproduct of 
its countable substructures. 
This explains why we can additionally restrict ourselves by 
countable structures.  

Theorem \ref{fingr} below shows that $\mathcal{G}_{w.sof}$ 
differs from the classes
$\{ \mathcal{G}_{sof} , \mathcal{G}_{hyplin},  \mathcal{G}_{l.sof} \}$.  
Is every  metric group a weakly sofic metric group? 
As above we may only consider countable discrete 
metric groups. 
The following remark shows that it suffices to consider the case 
of finitely generated free groups. 

\begin{remark} \label{doucha} 
{\em 
Concerning this question and Theorem \ref{main} M.Doucha 
has pointed out to the author that in fact any metric group 
embeds into a metric ultraproduct of finitely generated 
free groups with discrete bi-invariant metrics.  
This follows from the construction of \cite{Doucha} of a universal 
separable group $\mathbb{G}$ equipped with a complete 
bi-invariant metric bounded by 1. 
According to \cite{Doucha} $\mathbb{G}$ is the completion 
of a Fra\"{i}ss\'{e} limit of free metric groups as above. 
Thus the following question becomes principal in this direction. 
\begin{quote} 
Is every finitely generated free group with a bi-invariant metric 
weakly sofic? 
\end{quote} 
It is interesting that this is exactly a reformulation 
of Question 3.5 of \cite{Doucha}. 
Moreover it is also noted in \cite{Doucha} that this question is 
equivalent to extreme amenability of $\mathbb{G}$. }
\end{remark} 

\bigskip 

As we have already mentioned before there are examples 
which show that $\mathcal{G}_{sof}$ is a proper subset of 
$\mathcal{G}_{w.sof}$. 
Moreover these classes already have different finite members. 
It is the time for the presentation of the example. 

Let $p$ be a prime number $\ge 3$. 
Let us consider the cyclic group $\mathbb{Z}(p)$ with respect to 
so called Lee norm and Lee distances: 
$$ 
l_{Lee}(a) = \frac{2 min(a,p-a)}{p-1} \mbox{ , } d_{Lee}(a,b) = l_{Lee} (a-b) . 
$$ 
This is the normalized distance in the Cayley graph of $\mathbb{Z}(p)$. 
The latter was considered in \cite{ba}. 
Note that $l_{Lee} (1) = \frac{2}{p-1}$ is the minimal 
non-zero value and 
$$
l_{Lee} (\frac{p-1}{2}) = l_{Lee} (\frac{p+1}{2}) =1. 
$$ 

\begin{theorem} \label{fingr} 
Let $p\ge 13$. 
The metric group $(\mathbb{Z}(p), d_{Lee})$ 
does not belong to 
$\mathcal{G}_{sof} \cup \mathcal{G}_{hyplin} \cup \mathcal{G}_{l.sof}$. 
\end{theorem} 

{\em Proof.} 
Let us show that the group $(\mathbb{Z}(p), d_{Lee})$
is not $d_H$-approximable in the class of 
all symmetric groups $S_n$ with Hamming metrics $d_H$. 
Fix any generator $g$ of $\mathbb{Z}(p)$. 
In the situation of Definition \ref{DefMetrApp} 
let $\delta$ be sufficiently small, $n$ is much bigger 
than $p$ and let $F= \mathbb{Z}(p)\setminus \{ 0\}$. 
Note that for any function $\gamma$ from $F$ into $S_n$ 
as in that definition, any $\gamma (g^m )$ with $m\le p$ 
coincides with  $\gamma (g)^m$ on a $(1-p \cdot \delta )$-th
part of $\{ 1,..., n\}$ 
(apply the definition of $d_H$). 
Thus the function $\gamma$ defines an action of 
$\mathbb{Z}(p)$ on some $(1 - p\cdot \delta )$-th part 
of $\{ 1,...,n\}$. 
Since $\mathbb{Z}(p)$ is cyclic and simple, any non-trivial 
$f_1$ and $f_2\in F$ have the same orbits on this part 
of $\{ 1,...,n\}$. 
In particular $Fix(\gamma (f_1 )) = Fix(\gamma (f_2 ))$ 
there and the difference between Hamming norms of $\gamma (f_1 )$ 
and $\gamma (f_2 )$ does not exceed $p\cdot \delta$. 
Thus if $d_{Lee} (f_1 ,f_2 ) =\frac{1}{2}$  and 
$l_{Lee} (f_1 ) = \frac{2}{p-1}$, then the 
final inequality of Definition \ref{DefMetrApp} 
cannot hold for $f_1$ and $f_2$ until each 
$Fix (\gamma (f_1 ))$ and $Fix (\gamma (f_2 ))$  
is approximately a half of $\{ 1,...,n\}$. 
In the latter case that inequality does not hold for $0$ 
(the neutral element) and $f_1$. 

To see that $(\mathbb{Z}(p), d_{Lee})$ 
does not belong to $\mathcal{G}_{l.sof}$ 
we apply the same argument replacing 
$(S_n ,d_H )$ by $(GL_n (\mathbb{C}), \rho )$. 
If $g$, $\delta$ and $\gamma$ are as above, then 
we get that $\gamma (g)$ almost coincides 
with a semisimple element of $GL_n (\mathbb{C})$, 
say $a$ (because $\gamma (g)$ is $p\delta$-close 
with respect to $\rho$ to an element of finite order). 
We know that $a$ is presented by 
a diagonal matrix over $\mathbb{C}$. 
Since the diagonal elements of $a$ 
are $p$-th roots of $1$,  we see that 
when $0<i<p$ the matrices  
$a - id$ and $a^i - id$ have the same $0$-entries.  
This provides a contradiction as above. 

In the hyperlinear case we consider 
$(U(n), \frac{1}{2} d_{HS} )$ 
instead of $(GL_n (\mathbb{C}), \rho )$. 
Then fixing an appropriate base 
of $\mathbb{C}^n$ we find a diagonal 
matrix $a$ so that for all $0<j<p$ 
the powers $a^j$ are 
$p\delta$-close 
to $\gamma (g^j )$ with respect to 
$\frac{1}{2} d_{HS}$.  
Assume that
$l_{Lee} (g ) = \frac{2}{p-1}$ 
(otherwise we can choose another generator)  
and the $\frac{1}{2} d_{HS}$-distance between  
$a$ and $id$ is a number $\delta$-close  
to $\frac{2}{p-1}$. 
Note that the diagonal elements of $a$ 
are $p$-th roots of $1$, 
say $e^{i\phi_1} ,...,e^{i \phi_n}$.  
In particular each $a^j$ can be identified 
with the point 
$\frac{1}{\sqrt{n}}(e^{i\cdot j\phi_1} ,...,e^{i \cdot j\phi_n})$
of the unit ball in $\mathbb{C}^n$. 
Moreover the $d_{HS}$-distance 
between any $a^k$ and $a^l$ is exactly 
the distance between the corresponding points 
of the ball. 
Note that all $d_{HS} (a^j ,a^{j+1})$, 
$0\le j <p$, are the same. 
Since the module of 
the number $e^{i \phi_k }-1$ 
equals to $\sqrt{2 - 2 cos (\phi_k )}$, 
these distances can be computed as follows: 
$$ 
d_{HS} (a, id ) =
\frac{\sqrt{2}}{\sqrt{n}} \sqrt{1 - cos (\phi_1 )+ ...+ 1 - cos (\phi_n )} .
$$  
On the othe hand our assumptions imply 
$$ 
\frac{2}{p-1} \approx \frac{1}{2} d_{HS} (a, id ) .
$$  
Thus we see that $2$ 
($= \frac{4}{p-1} \cdot \frac{p-1}{2}$)  
can be approximately presented as 
$$
\sum^{(p-1)/2}_{j=1} d_{HS} (a^{j-1}, a^j) , 
$$ 
which in turn is the sum of distances between 
the corresponding points of the unit ball. 
On the other hand 
the $d_{HS}$-distance from $a^{(p-1)/2}$ to $id$ is 
$$  
\frac{\sqrt{2}}{\sqrt{n}} \sqrt{1 - cos ((p-1) \phi_1 /2)+ ...+ 1 - cos ((p-1) \phi_n /2)} .  
$$
Since $a^{(p-1)/2}$ does not represent 
$\frac{1}{\sqrt{n}}(-1 ,...,-1 )$,  
this distance is much shorter than the sum above 
(we may assume that the difference cannot be compared with 
$p\delta$).  
Thus $d_{Lee}$ cannot be approximated 
by $\frac{1}{2} d_{HS}$. 
$\Box$ 

\bigskip 

\begin{remark} 
{\em 
It is clear that the statement above holds 
if we replace $\mathbb{Z}(p)$ by any 
metric group where $(\mathbb{Z}(p), d_{Lee})$ 
can be embedded as a metric subgroup 
(with possible normalization). 
For example let us consider $\mathbb{Z}(p^s )$, 
where $p$ is prime and $s\ge 2$. 
Let $g$ be a generator of $\mathbb{Z}(p^s )$. 
Let us consider the sequence of embeddings 
$$
p^{s-1} \mathbb{Z}(p^s) \rightarrow 
p^{s-2} \mathbb{Z}(p^s )\rightarrow ... \rightarrow 
p \mathbb{Z}(p^s ) \rightarrow \mathbb{Z}(p^s ) . 
$$ 
If $g_1 , g_2 \in \mathbb{Z}(p^s )$ belong to distinct 
$p^k \mathbb{Z}(p^s )$-classes, where $k$ is minimal, 
let us define 
$$
d (g_1 ,g_2 )= \frac{1}{2^{k-1}}. 
$$  
If $g_1 - g_2 \in p^{s-1} \mathbb{Z}(p^s )$, then 
identifying the fixed generators of $p^{s-1}\mathbb{Z}(p^s )$ 
(i.e. $p^{s-1} g$) and of $\mathbb{Z}(p)$ 
we define the distance 
$$ 
d(g_1 ,g_2 ) = \frac{1}{2^{s-1}} l_{Lee} (g_1 - g_2 ) .
$$  
It is easy to see that $d$ is a metric, and 
behaves as an ultrametric for elements from distinct cosets. } 
\end{remark}

\section{Shifting metrics} 

In order to prove Theorem \ref{main} let us consider the following construction. 

\begin{definition} 
{\em 
Let $(G,d)$ be a metric group with $d\le 1$ and $\varepsilon \in [0,1]$. 
Let $d_{\varepsilon}$ be the metric on $G$ defined by the following rule: 
$$ 
d_{\varepsilon} (g,h) =  \frac{d(g,h)+\varepsilon }{1 + \varepsilon}  \mbox{ , for } g\not= h. 
$$ 
We call $d_{\varepsilon}$ the $\varepsilon$-shift of $d$.
\footnote{we are grateful to Krzysztof Majcher who suggested, 
that this construction is useful under our circumstances} 
}
\end{definition} 

It is easy to see that if $(G,d)$ is a bi-invariant metric group, 
then $(G, d_{\varepsilon})$ is a bi-invariant metric group too. 
Moreover $(G, d_{\varepsilon})$ is discrete. 
The hyperlinear case of Theorem \ref{main}(2) 
is based on the following theorem. 
We think that the sofic case below is interesting in itself.  

\begin{theorem} \label{shift} 
If $(G,d)$ is a sofic (resp.  hyperlinear) 
metric group  then $(G,d_{\varepsilon})$ is a sofic 
(resp. hyperlinear) metric group too. 
\end{theorem} 

{\em Proof.} 
Let us consider the case of 
a sofic metric group $(G,d)$. 
We will assume that $\varepsilon$ 
is sufficiently small. 
Take any  
$\delta < \frac{\varepsilon}{1+\varepsilon}$. 
For every finite $F\subset G$ 
there is a number $n$ and an embedding 
$\theta: F\rightarrow S_n$ so that  \\ 
- if $g,h,gh \in F$, then $d_{H} (\theta (g) \theta (h), \theta (gh)) <\delta$, \\ 
- if $1\in F$, then $d_H (\theta (1) , id ) < \delta$, \\  
- for all distinct $g,h \in F$, $|d(h,g) - d_{H}(\theta (h) , \theta (g) )| < \delta$.

On the other hand since $G$ is sofic by the proof of Theorem 3.5 of \cite{pestov} 
there are $m$ and an embedding $\theta' : F \rightarrow S_m$ so that \\  
- if $g,h,gh \in F$, then $d_{H} (\theta' (g) \theta' (h), \theta' (gh)) <\delta$, \\ 
- if $1\in F$, then $d_H (\theta' (1) , id ) < \delta$, \\  
- for all distinct $g,h \in F$, $|d_{H}(\theta' (h) , \theta' (g) )- 1|< \delta$. 

Note that $m$ can be taken arbitrary large: having $\theta'$ 
we can copy it as many times as we need. 
If $m'>m$ then extending elements of $F$ to 
the set $\{ m+1 ,m+2 ,...,m'\}$ by identity, 
we can reduce the numbers $\delta$ and $1 - \delta$  
in these conditions to a smaller ones. 
In particular taking $m'$ large enough considering 
representations in $S_{m'}$, 
we may modify the last contition  as follows: \\ 
- for all distinct $g,h \in F$, 
$|d_{H}(\theta' (h) , \theta' (g) )- \frac{\varepsilon}{1+\varepsilon}|< \frac{\delta \cdot \varepsilon}{1+\varepsilon}$. 

Note that this construction guarantees that 
$\frac{m}{m'}$ is close to  $\frac{\varepsilon}{1+\varepsilon}$. 
Since $m$ and $m'$ can be taken arbitrary large 
we may assume (neglecting small subsets) that $n|(m'-m)$ and 
moreover there is an equivalence relation $E$ on 
the set $\{ 1 ,...,m'\}$ consisting of classes of the same 
size so that $\{ 1,...,m\}$ and $\{ m+1 ,...,m'\}$ 
are $E$-invariant and $|\{ m+1 ,...,m'\} /E| =n$. 
Identifying $\{ 1,...,n\}$ with $\{ m+1 ,...,m'\}/E$ 
 we will assume that the action of elements of  
$F$ on $\{ m+1 ,...,m'\}$ is defined by $\frac{m'-m}{n}$ copies 
of the action $\theta$ of elements of $F$ on $\{ 1,...,n\}$ 
so that each copy is a transversal with respect to $E$.

This naturally amalgamates $\theta$ and 
$\theta' \upharpoonright \{ 1,...,m\}$ to 
permutations on $\{ 1,...,m' \}$.  
It is easy to see that the resulting map 
$\theta'' :F \rightarrow S_{m'}$ satisfies the following conditions: \\ 
- if $g,h,gh \in F$, then $d_{H} (\theta'' (g) \theta'' (h), \theta'' (gh)) <\delta$, \\ 
- if $1\in F$, then $d_H (\theta'' (e) , id ) < \delta$, \\  
- for all distinct $g,h \in F$, $|d_{\varepsilon}(h,g) - d_{H}(\theta'' (h) , \theta'' (g) )| < \delta$. 

To see the latter note that under the natural action of 
$h^{-1}g$ on $\{ m+1 ,...,m'\}$ the size of the support 
is approximately $d(g,h) \cdot (m'-m)$. 
Since $\frac{m'-m}{m'}$ is approximately 
$\frac{1}{1+\varepsilon }$ we see that 
$d_{H}(\theta'' (h) , \theta'' (g) )$ 
is approximately 
$$
\frac{d(g,h) }{1+\varepsilon } + \frac{\varepsilon }{1+\varepsilon } = 
\frac{d(g,h) + \varepsilon}{1+\varepsilon }. 
$$ 

We now discuss how this argument can be adapted to 
hyperlinear groups. 
Let us start with a metric hyperlinear group 
$(G,d)$ with $d\le 1$. 
Choosing $\varepsilon$, $\delta$ and a finite set $F$ 
as above we find a number $n$ and 
an embedding $\theta :F \rightarrow U(n)$ so that \\ 
- if $g,h,gh \in F$, then $d_{HS} (\theta (g) \theta (h), \theta (gh)) <\delta$, \\ 
- if $1\in F$, then $d_{HS} (\theta (1) , id ) < \delta$, \\  
- for all distinct $g,h \in F$, $|d(h,g) - \frac{1}{2}d_{HS}(\theta (h) , \theta (g) )| < \delta$. 

By Theorem 3.6 and Remark 3.7 of \cite{pestov} 
there are $m$ and an embedding 
$\theta' : F \rightarrow U(m)$ so that \\  
- if $g,h,gh \in F$, then $d_{HS} (\theta' (g) \theta' (h), \theta' (gh)) <\delta$, \\ 
- if $1\in F$, then $d_{HS} (\theta' (1) , id ) < \delta$, \\  
- for all distinct $g,h \in F$, $\frac{1}{2}d_{HS}(\theta' (h) , \theta' (g) )> 1- \delta$.  

To extend the remaining argument from the sofic 
to the hyperlinear case we need the following 
general observation. 
Given   
$\mu : F \rightarrow U(k)$,   
which embeds $F$ into $U(k)$ as above,  
one can naturally copy the actions of elements 
of $F$ in the space $\mathbb{C}^{k'k}$ 
so that the distance between the corresponding 
images of $g$ and $h$ in $U(k'k)$ 
coincides with 
$\frac{1}{2}d_{HS}(\mu (h) , \mu (g) )$. 
On the other hand if instead of this we extend 
$\mu$ to the remaining direct summands in 
$\mathbb{C}^{k'k}$ by identity 
we can reduce the distances as 
much as we need. 
Thus the argument above can be carried 
over in the hyperlinear case.  
$\Box$ 

\bigskip 

{\em Proof of  Theorem \ref{main}.} 
Let $\{ \varepsilon_n \} \rightarrow 0$. 
Let $\mathcal{U}$ be a non-principal ultrafilter on $\omega$. 
Then $(G,d)$ embeds into 
$\prod_{i\in \omega} (G,d_{\varepsilon_i})/\mathcal{U}$ 
under the diagonal map. 
Since $(G,d)$ is complete, it is 
a closed subgroup of the ultraproduct. 
The rest follows from Theorem \ref{shift}.  
$\Box$  

\bigskip 
 
In the linear sofic case  Theorem \ref{shift} 
should be modified. 
The problem is that  in this case the corresponding version of 
amplification does not work in the same way as in the case 
of sofic and hyperlinear groups. 
Namely it is not known if the value of the function $\alpha$ 
from Definition \ref{DefApp} can be taken 
equal to any real number strictly between $0$ and $1$.  
It is only proved in Proposition 5.13 of \cite{ArP2} 
that it can be made $\ge 1/4$:  
\begin{quote} 
a group $G$ is linearly sofic if and only if it can be 
embedded into a metric ultraproduct of metric 
groups $(GL_n (\mathbb{C}),\rho )$ so that 
any non-identity element of $G$ 
is distant from the identity by $\ge \frac{1}{4}$. 
\end{quote} 
We denote the metric in this ultraproduct by $d_{\omega}$ 
and consider $G$ under the corresponding restriction 
of $d_{\omega}$. 
Let us define the following metric on $G$:  
$$ 
d^{\omega}_{\varepsilon} (g,h) =  \frac{d(g,h)+\varepsilon d_{\omega}(g,h)}{1 + \varepsilon}  \mbox{ , when } g\not= h. 
$$ 
It is clear that for $g\not=h$,  
$$ 
d^{\omega}_{\varepsilon} (g,h) \ge  \frac{\varepsilon }{4 + 4\varepsilon} ;
$$ 
thus $(G, d^{\omega}_{\varepsilon})$ is a discrete group. 

\begin{theorem} \label{linshift} 
If $(G,d)$ is a linearly sofic  
metric group  then $(G,d^{\omega}_{\varepsilon})$ 
is a linearly sofic  metric group too. 
\end{theorem}

{\em Proof.} 
We adapt the proof of Theorem \ref{shift} as follows. 
Take any  
$\delta < \frac{\varepsilon}{1+\varepsilon}$. 
For every finite $F\subset G$ 
there is a number $n$ and an embedding 
$\theta: F\rightarrow GL_n (\mathbb{C})$ so that  \\ 
- if $g,h,gh \in F$, then $\rho (\theta (g) \theta (h), \theta (gh)) <\delta$, \\ 
- if $1\in F$, then $\rho (\theta (1) , id ) < \delta$, \\  
- for all distinct $g,h \in F$, $|d(h,g) - \rho (\theta (h) , \theta (g) )| < \delta$.   

On the other hand there are $m$ and 
an embedding $\theta' : F \rightarrow GL_m (\mathbb{C})$ so that \\  
- if $g,h,gh \in F$, then $\rho (\theta' (g) \theta' (h), \theta' (gh)) <\delta$, \\ 
- if $1\in F$, then $\rho (\theta' (1) , id ) < \delta$, \\  
- for all distinct $g,h \in F$, 
$| \rho (\theta' (h) , \theta' (g) )- d_{\omega} (h,g)|< \delta$. 

The number $m$ can be taken arbitrary large:  
for some $k$ one can naturally copy the actions 
of elements of $F$ in the space $\mathbb{C}^{mk}$ 
so that the distance between the corresponding 
images of $g$ and $h$ in $GL_{km}(\mathbb{C})$ 
coincides with 
$\rho(\theta' (h) , \theta' (g) )$. 

If $m'>m$ then extending  
$\theta$ to the remaining direct summands in 
$\mathbb{C}^{m'}$ by identity 
we can reduce the numbers $\delta$ and  $d_{\omega}(h,g) - \delta$  
in these conditions. 
We make the reduction so that: \\ 
- for all distinct $g,h \in F$, 
$|\rho (\theta' (h) , \theta' (g) )- \frac{\varepsilon d_{\omega}(h,g)}{1+\varepsilon}|< \frac{\delta\cdot \varepsilon}{1+\varepsilon}$. 

As in the proof of Theorem \ref{shift},   
$\frac{m}{m'}$ is close to  $\frac{\varepsilon}{1+\varepsilon}$. 
Applying the remaining argument of Theorem \ref{shift} 
we obtain a map 
$\theta'' :F \rightarrow GL_{m'}(\mathbb{C})$ so that 
for all distinct $g,h \in F$, 
$\rho (\theta'' (h) , \theta'' (g) )$ 
is approximately 
$$
\frac{d(g,h) }{1+\varepsilon } + \frac{\varepsilon d_{\omega}( g,h)}{1+\varepsilon }. 
$$ 
$\Box$ 

\bigskip

Note that Theorems \ref{main} and \ref{shift} 
obviously hold in the case of any variety of groups. 
For example we can everywhere replace hyperlinearity 
by nilpotence of degree $n$ or solubility of degree $n$. 

On the other hand it is worth noting that many classes 
of metric groups studied  in geometric group theory  
do not satisfy the  property of Theorem \ref{shift}. 
Indeed for example any compact group can be considered 
with respect to a bi-invariant metric $d \le 1$. 
Such a group is always amenable and has Kazhdan's property ${\bf (T)}$,  
\cite{BHV}. 
Since there are compact groups containing non-abelian free subgroups,  
moving from $d$ to $d_{\varepsilon}$ we can lose amenability.  
We can also consider uncountable compact groups $(G,d)$ 
(or just compact groups without finite Kazhdan's sets). 
In this case $(G,d_{\varepsilon})$ does not have ${\bf (T)}$.    
Adding an appropriate direct summand these examples can be made 
non-compact locally compact. 

Note that it is easy to see that if $(G,d_{\varepsilon})$ 
is amenable (has ${\bf (T)}$, ${\bf FH}$ or ${\bf OB}$, 
see \cite{BHV}, \cite{cornulier}), 
then so is $(G,d)$.  
On the other hand it is worth noting that all these properties 
and their negations are not axiomatizable in continuous logic 
(see \cite{sasza}). 

\begin{remark}  
{\em 
We do not know if any version of Theorem \ref{shift} 
holds in the case of weakly sofic metric groups. 
For example, we do not know if the condition that $(G,d)$ 
is weakly sofic implies that $(G, d_{\varepsilon})$ 
is weakly sofic too.  
At first sight the following idea may work.  
Let $(G,d)$ embed into 
a metric ultraproduct of finite metric
groups $(G_i ,d_i )$. 
Then one can try to approximate $(G, d_{\varepsilon})$
by all $(G_i , (d_i )_{\varepsilon})$. 
The problem is that elements of $\prod G_i$ 
of norm $0$ under metrics $d_i : i\in \mathbb{N}$,  
have norm $\varepsilon$ under metrics 
$(d_i )_{\varepsilon} : i\in \mathbb{N}$. 
Thus we do not even know if the metric ultraproduct 
of  all $(G_i , (d_i )_{\varepsilon})$ has 
a subgroup isomorphic to $G$. }
\end{remark}

\bigskip

1. Institute of Mathematics, University of Wroclaw, pl.Grunwaldzki 2/4, 

50-384 Wroclaw, Poland  

2. Institute of Mathematics, Silesian University of Technology, \parskip0pt

ul. Kaszubska 23, 48-101 Gliwice, \parskip0pt

Poland  \parskip0pt

{\em Aleksander.Iwanow@polsl.pl} 
\end{document}